\documentclass[12pt]{amsart}
\usepackage{amssymb,amsmath,epsfig,amscd,eucal}

\newtheorem{thm}{Theorem}

\newtheorem{lem}[thm]{Lemma}

\newtheorem{rem}[thm]{Remark}
\newtheorem{cor}[thm]{Corollary}
\newtheorem{exa}[thm]{Example}

\newtheorem{letterthm}{Theorem}

\newenvironment{pf}{\par\medskip\noindent{\em Proof. }}{\hfill $\square$\par\medskip}

\newcommand{\F}{\mathbb{F}}

\title{Subgroup separability in residually free groups}

\author[Martin R.~Bridson]{Martin R.~Bridson}
\address{Mathematics Department,
Imperial College London,
SW7 2AZ, UK}

\email{m.bridson@imperial.ac.uk}

\author[Henry Wilton]{Henry Wilton}
\address{Department of Mathematics, 1 University Station C1200, Austin, TX 78712-0257, USA}

\email{henry.wilton@math.utexas.edu}

\date{28th June 2007}

\begin{document}

\begin{abstract}
We prove that the finitely presentable subgroups of residually free groups are separable and that the subgroups of type $\mathrm{FP}_\infty$ are virtual retracts. We describe a uniform solution to the membership problem for finitely presentable subgroups of residually free groups.
\end{abstract}

\maketitle

\section{Introduction}

The importance of residual finiteness and subgroup separability have long been recognized, both in group theory and topology. The importance of residual freeness, on the other hand, has only come into focus recently despite being studied for several decades, cf.~\cite{Baum67}. This increased awareness of residually free groups is due to the central role that they play in the resolution of Tarski's problems on the logic of free groups (see \cite{Se1}, \cite{Se2} \emph{et seq.} and \cite{KM98a}, \cite{KM98b} \emph{et seq.}). In that body of work, the central objects of study are the finitely generated \emph{fully residually free} groups, also known as {\em limit groups}.

Fix a non-abelian free group $\F$.  A group $G$ is \emph{residually free} if, given an element $g\in G\smallsetminus 1$, there exists a homomorphism $f:G\to\F$ with $f(g)\neq 1$.
And $G$ is \emph{fully residually free} if, for every finite subset $X\subseteq G\smallsetminus 1$, there exists a homomorphism $f:G\to\F$ with $1\notin f(X)$.  Sela \cite{Se1} and Baumslag--Myasnikov--Remeslennikov \cite{BMR} show that every finitely generated residually free group is a subgroup of a direct product of finitely many fully residually free groups.

Free groups and free abelian groups are fully residually free. So too are  the fundamental groups of compact surfaces of Euler characteristic at most $-2$. A direct product of residually free groups is residually free but the direct product of a non-abelian group and a non-trivial group is never fully residually free. The reader may wish to keep direct products of free and surface groups in mind as examples of residually free groups.

The second author \cite{Wilt} established strict constraints on the manner in which finitely generated subgroups can sit inside limit groups (see Theorem \ref{Wilton theorem}), while the first author, in joint work with Howie, Miller and Short \cite{BHMS1}, proved that the finiteness properties of residually free groups govern the manner in which they can sit inside direct products of limit groups (see Theorems \ref{BHMS theorem} and \ref{Nilpotent quotients}). Following some reductions, the results of  \cite{BHMS1} allow one to treat an arbitrary finitely presented residually free group as if it were a subdirect product of limit groups that contains a term of the lower central series. Thus, provided that one has lemmas to cover the required reductions, one can prove interesting facts about residually free groups by pulling back results about nilpotent groups.  In the context that interests us here---the results of \cite{Wilt} concerning subgroup separability and virtual retractions---we shall see that the required reductions pose little problem.

\begin{letterthm}\label{Theorem A}
If $G$ is a finitely generated residually free group and $H\subset G$ is a subgroup of type $\mathrm{FP}_\infty$ over $\mathbb Q$ then $H$ is a virtual retract of $G$.
\end{letterthm}

\begin{letterthm}\label{Theorem B}
If $G$ is a finitely generated residually free group and $H$ is a finitely presentable subgroup of $G$ then $H$ is separable in $G$.
\end{letterthm}

Theorem \ref{Theorem B} is new even in the case where $G$ is a direct product of free groups or surface groups.

We remind the reader of the meaning and origin of the terms appearing in these theorems.  A subgroup $H\subset G$ is a \emph{virtual retract} if there exists a finite-index subgroup $V\subset G$, containing $H$, and a homomorphism $\rho : V\to H$ that restricts to the identity on $H$.  Following Long and Reid \cite{LR}, if every finitely generated subgroup of $G$ is a virtual retract then we say that $G$ has \emph{local retractions} or \emph{property LR}.

A finitely generated group $G$ is \emph{of type $\mathrm{FP}_n$ over $\mathbb Q$} if $\mathbb Q$, viewed as a trivial  $\mathbb Q G$-module,  has a projective resolution in which the first $(n+1)$ resolving modules are finitely generated. The group $G$ is \emph{of type $\mathrm{FP}_\infty$ over $\mathbb Q$} if it is  $\mathrm{FP}_n$ over $\mathbb Q$ for all $n$.  If $G$ has a compact $K(G,1)$ then it is of type $\mathrm{FP}_\infty$, and the converse holds for residually free groups \cite{BHMS1}.

A subgroup $H\subset G$ is termed \emph{separable} if, for every $g\in G\smallsetminus H$, there exists a finite-index subgroup $K\subset G$ containing $H$ but not $g$.
\begin{rem}\label{Pullbacks preserve separability}
If $f:G'\to G$ is an epimorphism of groups then $H\subset G$ is separable in $G$ if and only if $f^{-1}(H)$ is separable in $G'$.
\end{rem}
If every finitely generated subgroup of $G$ is separable then $G$ is called \emph{locally extended residually finite}, or \emph{LERF}.
An old theorem of Mal'cev shows that finitely generated
virtually nilpotent groups are LERF.


Already in the direct product of three free groups one finds finitely subgroups that are not of type $\mathrm{FP}_\infty$; the first examples are due to Stallings \cite{St}. Thus Theorem \ref{Theorem A} would not remain valid if we assumed only that $H$ was finitely presented. Likewise, the hypothesis of Theorem \ref{Theorem B} cannot be weakened to allow for finitely generated subgroups, because such subgroups are not separable in general, even in the direct product of two free groups, as the following example shows.

\begin{exa}
Let $Q=\langle a_1,\ldots,a_m|r_1,\ldots,r_n\rangle$ be a finite presentation for a non-residually finite group, for example
\[
BS(2,3)=\langle a,b|~a^{-1}b^2a=b^3\rangle,
\]
and let $q:F\to Q$ be the corresponding surjection from a free group $F=\langle a_1,\ldots,a_n\rangle$.  The diagonal subgroup $\Delta$ is not separable in $Q\times Q$.  Thus $H=(q\times q)^{-1}(\Delta)$ is not separable in $F\times F$, by Remark \ref{Pullbacks preserve separability}.  But $H$ is generated by the finite set of elements $\{(a_i,a_i)|1\leq i\leq m\} \cup\{(r_j,1)|1\leq j\leq n\}\subset F\times F$.
\end{exa}

Finally, we would like to point out that although it is far from view in the present article, the two trains of thought that we harness to prove our results (\cite{Wilt} and \cite{BHMS1}) both draw crucially on ideas that originate in Stallings' proof \cite{St2} of Marshall Hall's Theorem, i.e.\ Theorem \ref{Theorem A} for free groups.

\section{Reduction to subdirect products}

The purpose of this section is to reduce Theorems \ref{Theorem A} and \ref{Theorem B} to the special case where $G$ is a direct product of limit groups and $H\subset G$ is a subdirect product. (We remind the reader that a subgroup of a direct product $S\subset G_1\times\dots\times G_n$ is termed a {\em subdirect product} if its projection to each $G_i$ is onto.)

Our results rely crucially on the following result of Sela \cite[Claim 7.5]{Se1} and Baumslag--Myasnikov--Remeslennikov \cite[Corollary 19]{BMR}.

\begin{thm}\label{RF groups are subgroups of direct products}
A finitely generated group $G$ is residually free if and only if it is a subdirect  product of finitely generated, fully residually free groups.
\end{thm}

If $H\subset G$ is separable in (respectively, is a virtual retract of) $G$ and $K\subset G$ is any subgroup then $H\cap K$ is separable in (respectively, is a virtual retract of) $K$.  It follows that we need only consider the case when
\[
G=\Lambda_1\times\cdots\times \Lambda_n
\]
where each $\Lambda_i$ is a limit group.  Write $\pi_i:G\to \Lambda_i$ for projection onto the $i$th factor.

In order to apply the results of \cite{BHMS1} we must also arrange that the subgroups $H\subset G$ that we consider are subdirect products. This is achieved by using the following lemma to replace $G$ by $\hat{H}:=\pi_1(H)\times\cdots\times \pi_n(H)$. The proof of this lemma applies to subdirect products of arbitrary groups with property LR. Its applicability in our situation rests on work of the first author.

\begin{thm}[Theorems A \& B of \cite{Wilt}]\label{Wilton theorem}
If $G$ is a finitely generated, fully residually free group and $H$ is a finitely generated subgroup of $G$ then $H$ is separable in $G$ and a virtual retract of $G$.
\end{thm}

\begin{lem}\label{Reduction to subdirect products}
A subgroup $H$ is separable in (respectively, is a virtual retract of) $G$ if and only if it is separable in (respectively, is a virtual retract of) $\hat{H}$.
\end{lem}
\begin{pf}
It is easy to see that, if $H$ is separable in (respectively, is a virtual retract of) $G$ then $H$ is separable in (respectively, is a virtual retract of) $\hat{H}$.  So we will concentrate on the converse assertions.

For each $i$, $\pi_i(H)$ is a virtual retract of $\Lambda_i$, so there exists a finite-index subgroup $V_i\subset\Lambda_i$ and a retraction $\rho_i:V_i\to\pi_i(H)$.  The product retraction
\[
\rho=\rho_1\times\cdots\times\rho_n:V_1\times\cdots\times V_n\to \hat{H}.
\]
exhibits $\hat{H}$ as a virtual retract.  It follows immediately that if $H$ is a virtual retract of $\hat{H}$ then $H$ is a virtual retract of $G$.  Suppose $H$ is separable in $\hat{H}$ and $g\notin H$.  There are two cases to consider.  If $g\notin\hat{H}$ then, for some $i$, $\pi_i(g)\notin\pi_i(H)$.  As $\pi_i(H)$ is separable in $\Lambda_i$ it follows that $g$ can be separated from $H$.  If $g\in\hat{H}$ then, by the hypothesis that $H$ is separable in $\hat{H}$, there exists a finite-index subgroup $\hat{V}\subset\hat{H}$ containing $H$ but not $g$.  Now $\rho^{-1}(\hat{V})$ is a finite-index subgroup of $G$ containing $H$ but not $g$.
\end{pf}

With this lemma in hand, we shall henceforth
assume that $G=\Lambda_1\times\cdots\times \Lambda_n = \hat H$.

The following notation will be useful in what follows.  Let $A_i=\Lambda_1\times\cdots\times \Lambda_i$, let $B_i=\Lambda_{i+1}\times\cdots\times \Lambda_n$, and write
\[
p_i:G=A_i\times B_i\to A_i
\]
for the projection that kills $B_i$ and
\[
q_i:G=A_i\times B_i\to B_i
\]
for the projection that kills $A_i$.

\section{Virtual retractions}

The reductions of the previous section allow us to appeal to the following theorem of Bridson, Howie, Miller and Short \cite{BHMS1}.

\begin{thm}\label{BHMS theorem}
Let $D=\Lambda_1\times\cdots\times \Lambda_n$ be a direct product of non-abelian limit groups.  If a subdirect product $S\subset D$ is of type $\mathrm{FP}_n$ over $\mathbb Q$ and $\Lambda_i\cap S\neq 1$ for $i=1,\dots,n$, then $S$ has finite index in $D$.
\end{thm}

The following lemma completes the proof of Theorem \ref{Theorem A}.

\begin{lem}\label{Virtual retraction}
Let $G=\Lambda_1\times\cdots\times \Lambda_n$ be a product of limit groups and let $H\subset G$ be a subgroup that is of type $\mathrm{FP}_n$ over $\mathbb Q$.  Then $H$ is a virtual retract of $G$.
\end{lem}
\begin{pf}
According to Lemma \ref{Reduction to subdirect products}, we may assume that $H$ is a subdirect product.  We reorder the factors so that, for some $i$, $A_{i-1}\cap H$ is trivial and, for all $j\geq i$, $H\cap\Lambda_j\neq 1$. There is no loss of generality in assuming that for at most one $j\geq i$, say $j=i$, $\Lambda_j$ is abelian.  Decompose $\Lambda_j$ as $P\times Q$ where $H\cap P=1$ and $H\cap Q$ is of finite index in $Q$.  We relabel the factors again so that $P$ is subsumed into $A_{i-1}$ and $Q$ becomes $\Lambda_i$.

We have now arranged that $G=\Lambda_1\times\cdots\times\Lambda_n$ where, for some $i$, $A_{i-1}\cap H=1$, $\Lambda_i$ is abelian
(possibly trivial) with $H\cap\Lambda_i$ of finite index in $\Lambda_i$, and $B_i$ is a direct product of non-abelian limit groups that intersect $H$ non-trivially.  The projection $q_{i-1}$ is injective on $H$; let $\iota:q_{i-1}(H)\to G$ be the inverse of $q_{i-1}|_H$.

We claim that the intersection $S:=q_{i-1}(H)\cap B_i$ is finitely generated. Indeed $I:=q_{i-1}(H)\subset\Lambda_i\times B_i$ is finitely generated and intersects $\Lambda_i$ in a subgroup of finite index, so  $I_0:=(I\cap \Lambda_i) \times S$, being of finite index in $I$, is also finitely generated, and $S$ is a quotient of $I_0$.

Now $S$ satisfies the hypotheses of Theorem \ref{BHMS theorem} in the direct product $B_i=\Lambda_{i+1}\times\cdots\times \Lambda_n$, and hence is of finite-index. Therefore $q_{i-1}(H)$ has finite index in $\Lambda_i\times B_i= B_{i-1}$ and the restriction of $\iota\circ q_{i-1}$ defines a retraction $A_{i-1}\times q_{i-1}(H)\to H$.
\end{pf}

\section{Separability}

The proof that finitely presented subgroups are separable requires a further result from \cite{BHMS1}.

\begin{thm}[cf.~Theorem 4.2 of \cite{BHMS1}]\label{Nilpotent quotients}
Let $D=\Lambda_1\times\cdots\times\Lambda_n$ be a direct product of limit groups and suppose that $S\subset D$ is a finitely presented subdirect product with $S\cap\Lambda_i\neq 1$ for each $i$.  Then $S\cap\Lambda_i\lhd\Lambda_i$ and $\Lambda_i/(S\cap\Lambda_i)$ is virtually nilpotent for each $i$.
\end{thm}

That $S\cap\Lambda_i$ is normal in $\Lambda_i$ follows directly from the fact that $S$ is a subdirect product. In the argument of \cite{BHMS1},  Theorem \ref{Nilpotent quotients} appears under some additional hypotheses: that $n\geq 2$; that each $\Lambda_i$ is non-abelian; and that each $\Lambda_i$ splits as an HNN-extension over a cyclic subgroup, with stable letter $t_i\in S\cap\Lambda_i$.  But in the way we have stated the theorem,
each of these hypotheses can be removed.  If $n=1$ then, because $S$ is a subdirect product, $S=\Lambda_i$ so the quotient is trivial.  If $\Lambda_i$ is abelian then $\Lambda_i/(S\cap\Lambda_i)$ is abelian and therefore nilpotent.  As observed in the proof of part 5 of Proposition 3.1 of \cite{BHMS1}, each $\Lambda_i$ contains a finite-index subgroup $\Delta_i$ that decomposes as an HNN-extension of the required form.  As $\Delta_i/(S\cap\Delta_i)$ is virtually nilpotent, $\Lambda_i/(S\cap\Lambda_i)$ is also virtually nilpotent.

\begin{cor}\label{Nontrivial intersections are separable}
If $D$ and $S$ satisfy the hypotheses of Theorem \ref{Nilpotent quotients} then $S$ is separable in $D$.
\end{cor}
\begin{pf}
Let $N=(S\cap\Lambda_1)\times\cdots\times (S\cap\Lambda_n)$, a normal subgroup of $D$.  By Theorem \ref{Nilpotent quotients}, $D/N$ is virtually nilpotent.  Let $q:D\to D/N$ be the quotient map.  As virtually nilpotent groups are subgroup separable, $q(S)$ is separable in $D/N$.  Since $N\subset S$ we have $S=q^{-1}(q(S))$, so the result follows by Remark \ref{Pullbacks preserve separability}.
\end{pf}

In the light of Lemma \ref{Reduction to subdirect products}, the following lemma completes the proof of Theorem \ref{Theorem B}.

\begin{lem}\label{Separable}
Let $G=\Lambda_1\times\cdots\times\Lambda_n$ be a direct product of limit groups and let $H\subset G$ be a finitely presentable subdirect product.  Then $H$ is separable in $G$.
\end{lem}
\begin{pf}
Renumbering the factors, as before, we may assume there
 exists an integer $i$ so that $q_i|_H$ is a monomorphism and  $q_i(H)\cap\Lambda_j\neq 1$ for each $j>i$.  Let $\iota:q_i(H)\to H$ be the inverse of $q_i|_H$.  Let $g\in G\smallsetminus H$. We aim to find a finite-index subgroup of $G$ containing $H$ but not $g$.

Suppose that $g\notin A_i\times q_i(H)$.  As $q_i(H)\subset B_i$ satisfies the hypotheses of Corollary \ref{Nontrivial intersections are separable},  there exists a finite-index subgroup $K\subset B_i$ containing $q_i(H)$ but not $q_i(g)$, so $q_i^{-1}(K)$ separates $g$ from $H$ in $G$.

Suppose now that $g\in A_i\times q_i(H)$.
We first  claim that there exists a finite-index subgroup $V$ of $A_i\times q_i(H)$ containing $H$ but not $g$.  To see this, we write $H$ as the graph of $\iota$:
\[
H=\{(a,b)\in A_i\times q_i(H)|~a=\iota(b)\}.
\]
Consider the map $f:A_i\times q_i(H)\to A_i\times A_i$ defined by $(a,b)\mapsto (a,\iota(b))$.  The pre-image of the diagonal subgroup of $A_i\times A_i$ under $f$ is precisely $H$---in particular, $f(g)\notin f(H)$.  The claim now follows because diagonal subgroups of direct products of residually finite groups are always separable.  For if $f(g)=(\alpha_1,\alpha_2)$ with $\alpha_1\neq\alpha_2$, then since $A_i$ is residually finite there exists a map $\phi$ from $A_i$ to a finite group $Q$ such that $\phi(\alpha_1)\neq\phi(\alpha_2)$, so
\[
\phi\times\phi:A_i\times A_i\to Q\times Q
\]
is a map to a finite group and $(\phi\times\phi)\circ f(g)\notin (\phi\times\phi)\circ f(H)$.  If $V$ is the pre-image of the diagonal subgroup of $Q\times Q$ under the map $(\phi\times\phi)\circ f$, then $V$ is a finite-index subgroup of $A_i\times q_i(H)$ containing $H$ but not $g$.

It now suffices to separate $V$ from $g$ in $G$.  Notice that $V$ is a finitely presented subgroup of $G$ and, for each $i$, $V\cap\Lambda_i\neq 1$.  Furthermore, since $H$ is a subdirect product of $G$ and $H\subset V$, we have that $V$ is also a subdirect product.  Therefore $V$ satisfies the hypotheses of Theorem \ref{Nilpotent quotients}, so is separable.
\end{pf}

\section{The membership problem}

A solution to the membership problem for an arbitrary finitely presentable
subgroup of a product of free and
surface groups is given in section 5.2 of \cite{BM}.
The results of \cite{BHMS1} allow one to extend this to
finitely presented subgroups of residually free groups \cite{BHMS2}.
However, that solution is not uniform---the algorithm depends
on additional information about the subgroup. Here we provide a uniform solution.

\begin{cor}
There is an algorithm that, given a finite presentation for a residually free group $G\cong \langle A\mid R\rangle$, a finite set of words $S$ in the alphabet $A\cup A^{-1}$ and an additional word $g$, will determine whether or not $g$ lies in the subgroup of $G$ generated by $S$ provided that $S$ is finitely presentable. (If $S$ is not finitely presentable, the algorithm will terminate with the correct conclusion if $g\in\langle S\rangle$ but may fail to terminate if $g\notin\langle S\rangle$.)
\end{cor}

The proof is a standard argument about separable subgroups, which we outline here.  The algorithm consists of two processes run in parallel, one of which will reach a conclusion: on the one hand, one enumerates the finite quotients of $G$ and checks to see if the image of $g$ lies in the image of $\langle S\rangle$, stopping if it does not and declaring that $g\notin \langle S\rangle$; on the other hand, working in the free group $F(A)$, one proceeds along an enumeration $u_n$ of the words in the free monoid on $S\cup S^{-1}$ and a na\"ive enumeration of the products $p_m$ in $F(A)$ of conjugates of relators $r\in R$, checking to see if each $g^{-1}u_n$ is freely equal to $p_m$. (One proceeds diagonally through the enumerations $u_n$ and $p_m$.)

A limitation of the utility of this algorithm comes from the observation that in general there is no way of knowing if $\langle S\rangle$ is finitely presentable.  This is the case when $G$ is the direct product of two free groups, for example.

\bibliographystyle{plain}

\end{document}